\documentclass[12pt]{amsart}
\usepackage{amssymb, amscd, amsmath, amsthm, epsf, epsfig, latexsym, enumerate}

\newtheorem{theorem}{Theorem}
\newtheorem{lemma}[theorem]{Lemma}
\newtheorem*{cor}{Corollary}

\begin{document}

\title{$2$-Knots with solvable group}

\author{Jonathan A. Hillman }
\address{School of Mathematics and Statistics F07\\
     University of Sydney, Sydney\\ 
     NSW 2006, Australia }

\email{jonh@maths.usyd.edu.au}

\begin{abstract}
We complete the TOP classification of 2-knots with torsion-free, 
solvable knot group by showing that fibred 2-knots
with closed fibre the Hantzsche-Wendt flat 3-manifold $HW$
are not reflexive, while every fibred 2-knot with closed fibre
a $Nil$-manifold with base orbifold $S(3,3,3)$ is reflexive,
and by giving explicit normal forms for the strict weight orbits 
of normal generators for the groups of all knots in either class.
We also determine when the knots are amphicheiral or invertible,
and show that the only doubly null-concordant knots 
with such groups are those with group $\pi\tau_29_{46}$.
\end{abstract}



\maketitle
The largest class of groups $\pi$ over which TOP surgery techniques 
in dimension 4 are known to hold is the class $SA$ obtained 
from groups of subexponential growth by extensions and increasing unions.
No such group has a noncyclic free subgroup.
The known 2-knot groups in this class are either torsion-free 
and solvable or have finite commutator subgroup.
(It seems plausible that there may be no others.
See Theorem 15.13 of \cite{Hi} and \S4 below 
for evidence in this direction.) 

If the group of a nontrivial 2-knot $K$ 
is torsion-free and elementary amenable then $K$ is either  
the Fox knot (Example 10 of \cite{Fo}) or is fibred,
with closed fibre $\mathbb{R}^3/\mathbb{Z}^3$,
the Hantzsche-Wendt flat 3-manifold $HW=\mathbb{R}^3/G_6$
or a $Nil^3$-manifold.
(See Theorem 2.)
Each such knot is determined up to Gluck reconstruction,
TOP isotopy and change of orientations by its group $\pi$ and weight orbit
(the orbit of a weight element under the action of $Aut(\pi)$).
This orbit is unique for the Fox knot and for the fibred knots 
with closed fibre $\mathbb{R}^3/\mathbb{Z}^3$
(the Cappell-Shaneson knots)  
or a coset space of the Lie group $Nil$.
In each of these cases the questions of amphicheirality,
invertibility and reflexivity have been decided.
(See \cite{Hi,Hi09,HP}.)

There are just two possible knot groups $G(\pm)$ 
realized by knots with closed fibre $HW$. 
No such knot is reflexive.
The 3-twist spin of the figure eight knot $\tau_34_1$ 
and its Gluck reconstruction $\tau_34_1^*$ (with group $G(+)$) are 
strongly $\pm$amphicheiral but not strongly invertible.
The remaining knots have closed fibre the
2-fold branched cover of $S^3$, branched over a Montesinos knot
$k(e,\eta)=K(0|e;(3,\eta),(3,1),(3,1))$, with $e$ even and $\eta=\pm1$.
These are all reflexive.
This class includes the 2-twist spins $\tau_2k(e,\eta)$
of such Montesinos knots,
which are strongly $+$amphicheiral but not invertible.
With the exception of four other knots with group $G(-)$
none of the other knots is amphicheiral or invertible.
In all cases we give explicit normal forms for the strict weight orbits.
This completes the TOP classification of such knots.
(However, it is not known whether they are all smooth
knots in the standard $S^4$.)

When the commutator subgroup of a 2-knot group
is finite the list of possible groups 
and weight orbits is known, 
but the surgery obstruction groups are large, 
and there are in general infinitely many TOP locally flat knots
with a given such group.
Thus it is reasonable to restrict attention to those which are fibred.
The closed fibre is then a spherical manifold $S^3/\pi'$.
In this case the question of reflexivity has been settled 
for 10 of the 17 possible families of such knots \cite{PS}.
It is likely that none of the remaining knots are reflexive, 
but this has not yet been confirmed.

In the final section we show that a nontrivial 2-knot $K$
such that $\pi{K}$ is torsion-free and solvable 
is TOP doubly slice if and only if $\pi{K}\cong\pi\tau_29_{46}$.
We observe also that if $K$ is 2-knot such that 
$(\pi{K})'$ is finite then it is TOP doubly slice if and only if
$\pi{K}\cong\pi\tau_53_1\cong{Z\times{I^*}}$.

\section{knot groups}

If $G$ is a group then $G'$ and $\zeta{G}$ shall denote 
the commutator subgroup and the centre of $G$, 
while $C_G(H)$ and $N_G(H)$ shall genote the centralizer 
and normalizer of the subgroup $H\leq{G}$, respectively. 

An automorphism $\phi$ of a group $G$ is {\it meridianal\/} if 
$\langle\langle{g^{-1}}\phi(g)\mid{g}\in{G}\rangle\rangle_G=G$.
When $G$ is finitely generated and solvable this holds if and only if 
$H_1(\phi)-1$ is an automorphism of the abelianization $H_1(G)$.
If $\phi$ and $\psi$ are meridianal automorphisms of $G$ then the semidirect
products $G\rtimes_\phi{Z}$ and $G\rtimes_\psi{Z}$ are isomorphic if and
only if the outer automorphism class $[\phi]$ is conjugate to $[\psi]$
or $[\psi]^{-1}$ in $Out(G)$.
There is an isomorphism preserving the stable letters of the HNN extensions if
and only if $\phi$ and $\psi$ are conjugate in $Aut(G)$. 
(See Lemma 1.1 of \cite{Hi}.) 

Let $t\in\pi=G\rtimes_\phi{Z}$ be an element whose
normal closure $\langle\langle{t}\rangle\rangle_\pi$ is the whole group.
Every such ``weight element" $w$ is of the form $w=gt$ or $w=(gt)^{-1}$, 
for some $g\in{G}$.
The strict weight orbit of $w$ is the set
$\{\alpha(w)\mid \alpha\in Aut(\pi),~\alpha(w)\equiv w~mod~G\}$.

If $\pi$ has a weight element then $G=\pi'$.
Let $c_x\in{Aut(G)}$ be the automorphism induced by conjugation 
by $x$ in $\pi$.
Two weight elements $t$ and $gt$ 
with $g\in{G}$ are in the same strict weight orbit 
if and only if there is an automorphism $\psi$ of $G$ 
such that $c_g=\psi.c_t.\psi^{-1}c_t^{-1}$,
by Theorem 14.1(3) of \cite{Hi}.
In particular, $c_t$ and $c_{gt}$ have the same order,
and $[\psi][\phi]=[\phi][\psi]$ in $Out(G)$.

If $\pi$ is a solvable $n$-knot group then $\pi''$ acts 
transitively on the set of weight elements representing 
a given generator $t$ of $\pi/\pi'$, 
by Theorem 14.1(4) of \cite{Hi}.
Unfortunately, this action does not usually induce an action 
on the set of weight orbits.
This result does not extend to all virtually solvable groups.
In particular, it does not hold for $Z\times{I^*}$.

Although it is possible to study the automorphism groups considered below
by purely algebraic means, 
we shall use embeddings in the appropriate affine groups
to guide the construction of homeomorphisms and isotopies. 

\section{self-homeomorphisms of knot exteriors}

We assume that the spheres $S^n$ are oriented.
Let $K:S^2\to{S^4}$ be a 2-knot with exterior $X$, and fix a homeomorphism
$\partial{X}\cong{S^2\times{S^1}}$ which is compatible with 
the orientations of the spheres $S^1$, $S^2$ and $X\subset{S^4}$.
Let $\tau(x,y)=(\rho(y)(x),y)$ for all $(x,y)$ in $S^2\times{S^1}$, 
where $\rho:S^1\to SO(3)$ is an essential map.
The {\it Gluck reconstruction\/} of $K$ is the knot $K^*$
given by the composite inclusion
\[S^2\subset{S^2\times{D^2}}\subset{X\cup_\tau{S^2\times{D^2}}}\cong{S^4}.\]
The knot $K$ is {\it reflexive\/} if $K^*$ is isotopic to one 
of the four knots $K,r_4K,Kr_2$ or $-K=r_4Kr_2$ obtained 
by composition with reflections $r_n$ of $S^n$.

We may extend any self-homeomorphism $h$ of $X$ ``radially" 
to a self-homeomorphism of the knot manifold
$M(K)=X\cup{D^3\times{S^1}}$ 
which maps the cocore $C=\{0\}\times{S^1}$ to itself.
If $h$ preserves both orientations or reverses both orientations 
then it fixes the meridian, and we may assume that $h|_C=id_C$.
If $h$ reverses the meridian $t$, 
we may still assume that it fixes a point on $C$.
We take such a fixed point as the basepoint for $M(K)$.
Let $h'_*$ be the induced automorphism of $\pi'$.

If $K$ is invertible or $\pm$amphicheiral there is a
self-homeomorphism $h$ of $(S^4,K)$ which changes the orientations
appropriately, but does not twist the normal bundle of $K(S^2)\subset{S^4}$. 
If it is reflexive there is such a self-homeomorphism which changes
the framing of the normal bundle.
Thus if $K$ is $-$amphicheiral there is such an $h$ which 
reverses the orientation of $M(K)$ and $h'_*$ commutes with the
meridianal automorphism $c_t$.
If $K$ is invertible or $+$ampicheiral there is a homeomorphism $h$
such that $h'_*c_th'_*=c_t^{-1}$ and which preserves or reverses the
orientation.

\section{sections of the mapping torus}

Let $\theta$ be a self-homeomorphism of a 3-manifold $F$,
with mapping torus $M(\theta)=F\times[0,1]/\thicksim$,
where $(f,0)\thicksim(\theta(f),1)$ for all $f\in{F}$,
and canonical projection $p_\theta:M(\theta)\to{S^1}$,
given by $p_\theta([f,s])=e^{2\pi{is}}$ for all $[f,s]\in{M(\theta)}$. 
The mapping torus $M(\theta)$ is orientable if and only if $\theta$
is orientation-preserving.
If $\theta'=h\theta{h^{-1}}$ for some self-homeomorphism $h$ of $F$
then $[f,s]\mapsto[h(f),s]$ defines a homeomorphism 
$m(h):M(\theta)\to{M(\theta')}$ such that $p_{\theta'}m(h)=p_\theta$.
Similarly, if $\theta'$ is isotopic to $\theta$ then
$M(\theta')\cong{M(\theta)}$.

If $P\in{F}$ is fixed by $\theta$ then the image of $P\times[0,1]$ in
$M(\theta)$ is a section of $p_\theta$.
In particular, if the fixed point set of $\theta$ is connected there is
a canonical isotopy class of sections.
If moreover $\theta_*=\pi_1(\theta)$ is meridianal these determine a preferred
conjugacy class of weight elements in the group $\pi_1(M(\theta))$.
(Two sections are isotopic if and only if they 
represent conjugate elements of $\pi$.)

In general, we may isotope $\theta$ to have a fixed point $P$.
Let $t\in\pi_1(M(\theta))$ correspond to the constant section of
$M(\theta)$, and let $u=gt$ with $g\in\pi_1(F)$.
Let $\gamma:[0,1]\to{F}$ be a loop representing $g$.
There is an isotopy $h_s$ from $h_0=id_F$ to $h=h_1$ which drags $P$
around $\gamma$, so that $h_s(P)=\gamma(s)$ for all $0\leq{s}\leq1$.
Then $H([f,s])=[(h_s)^{-1}(f),s]$ defines a homeomorphism
$M(\theta)\cong{M(h^{-1}\theta)}$.
Under this homeomorphism the constant section of
$p_{h^{-1}\theta}$ corresponds to the section of $p_\theta$
given by $m_u(t)=[\gamma(t),t]$, which represents $u$.
If $F$ is a geometric 3-manifold we may assume that $\gamma$
is a geodesic path.

Suppose henceforth that $\theta$ is orientation-preserving 
and $\theta_*$ is meridianal. 
Then surgery on a section gives a 2-knot.
There are two possible framings for the surgery, 
but the exteriors of the two knots are homeomorphic.

This is the situation for twist-spins,
where $F$ is a cyclic branched cover of $S^3$, 
branched over a classical knot,
and $\theta$ generates the covering group.
The subset fixed by $\theta$ is connected and nonempty,
since it is the branch locus.
The knot exterior is the complement of an open regular neighbourhood of the
canonical section of the mapping torus of $\theta$.

If $F$ has universal cover $\widetilde{F}\cong\mathbb{R}^3$
and $h$ is a self-homeomorphism of $M(\theta)$ which fixes 
a section setwise the behaviour of $h$ with respect to the orientations
is detected by the effect of $h'_*$ on $H_3(F;\mathbb{Z})$ and
whether $h'_*c_th'_*=c_t$ or $c_t^{-1}$.
As in \cite{CS,HP} (and Chapter 18 of \cite{Hi}), 
in order to determine whether $h$ changes the framing
it shall suffice to pass to the irregular covering space 
$M(\widetilde\theta)=\widetilde{F}\times_{\widetilde\theta}S^1$.
We seek a coordinate homeomorphism $\widetilde{F}\cong\mathbb{R}^3$ 
which gives convenient representations of the maps in question, 
and then use an isotopy from the identity to
$\widetilde\theta$ to identify $M(\widetilde\theta)$
with $\mathbb{R}^3\times{S^1}$.

\begin{lemma}
Let $K$ be a fibred $2$-knot. 
If there is a self homeomorphism $h$ of $X(K)$
which is the identity on one fibre and such that $h|_{\partial{X}}=\tau$ 
then all knots $\tilde{K}$ with $M(\tilde{K})\cong{M(K)}$ are reflexive.
In particular, this is so if the monodromy of $K$ has order $2$. 
\end{lemma}

\begin{proof}
We may extend $h$ to a self-homeomorphism $\hat{h}$ of $M(K)$
which fixes the surgery cocore $C\cong{S^1}$.
After an isotopy of $h$, we may assume that it is the identity on a product
neighbourhood $N=\hat{F}\times[-\epsilon,\epsilon]$ of the closed fibre.
Since any weight element for $\pi$ may be represented by a section
$\gamma$ of the bundle which coincides with $C$ outside $N$,
we may use $h$ to change the framing of the normal bundle 
of $\gamma$ for any such knot.
Hence every such knot is reflexive.

If the monodromy of $K$ has order 2 then
``turning the mapping torus upside-down" changes the 
framing of the normal bundle and fixes one fibre.
\end{proof}

The reflexivity of 2-twist spins is due to Litherland. See \cite{Li,Mo,Pl}. 

\section{torsion-free elementary amenable implies solvable}

We shall let $\Phi$ denote the group of the Fox knot.
This is an ascending HNN extension $\Phi\cong{Z*_2}$,
with presentation $\langle {a,t}\mid {tat^{-1}=a^2}\rangle$.

\begin{theorem}
Let $K$ be a $2$-knot whose group $\pi=\pi{K}$ is torsion-free 
and elementary amenable.
Then $K$ is trivial, the Fox knot, or is fibred with closed fibre
a flat $3$-manifold or a $Nil$-manifold.
\end{theorem}

\begin{proof}
If $\pi$ is torsion-free and has more than one end then $\pi\cong{Z}$,
and so $K$ is trivial \cite{Fr}.
If $\pi$ has one end and $H^2(\pi;\mathbb{Z}[\pi)=0$ 
then $M(K)$ is aspherical, by Theorem 3.5 of \cite{Hi}, 
and so $H^4(\pi;\mathbb{Z}[\pi)\not=0$.
Otherwise $H^2(\pi;\mathbb{Z}[\pi)\not=0$.
In all cases $H^s(\pi;\mathbb{Z}[\pi])\not=0$ for some $s\leq4$,
and so $\pi$ is virtually solvable, by Proposition 3 of \cite{Kr}.
It then follows that either $\pi\cong{Z}$ or $\pi\cong\Phi=Z*_2$, 
or that $\pi$ is virtually poly-$Z$ of Hirsch length 4.
(See Theorem 15.13 of \cite{Hi}.)
If $\pi\cong\Phi$ then $K$ is the Fox knot or its reflection
\cite{Hi09}, while the remaining cases are covered in Chapter 16 of \cite{Hi}.
\end{proof}

Can we relax the condition on torsion?
Let $\pi$ be an elementary amenable knot group. 
Since $\pi$ is finitely presentable and has an infinite cyclic quotient it is
an HNN extension with finitely generated base and associated subgroups.
Since it has no noncyclic free subgroups the HNN extension is ascending: 
$\pi\cong{H}*_\phi$, 
where $H$ is finitely generated and $\phi:H\to{H}$ is injective.
If moreover $H$ is $FP_3$ and virtually indicable then either $\pi'$ is finite
or $\pi$ is torsion-free, by Theorem 15.13 of \cite{Hi}.

The additional hypotheses on $H$ could be removed if we had a better
understanding of when $H^2(\pi;\mathbb{Z}[\pi])=0$.
Suppose that whenever $G$ is a finitely presentable group such that either
\begin{enumerate}
\item $G$ has an elementary amenable normal subgroup $E$ such that
\begin{enumerate}[(a)]
\item $h(E)>2$; or 
\item $h(E)=2$ and $G/E$ is infinite; or 
\item $h(E)=1$ and $G/E$ has one end; or
\end{enumerate}
\item $G\cong{B}*_\phi$ is an ascending HNN extension  
with finitely generated, 1-ended base $B$;
\end{enumerate}
then $H^2(G;\mathbb{Z}[G])=0$.

We may then argue as follows.
Since $\pi$ is finitely presentable and has an infinite cyclic quotient it is
an HNN extension with finitely generated base and associated subgroups.
Since it has no noncyclic free subgroups the HNN extension is ascending: 
$\pi\cong{H}*_\phi$, 
where $H$ is finitely generated and $\phi:H\to{H}$ is injective.
Since $\pi$ is elementary amenable and infinite $\beta_1^{(2)}(\pi)=0$.
If $h(\pi)=1$ then $\pi'$ is finite.
Suppose that $\pi'$ is infinite. 
Then $\pi$ has one end.
If $h(\pi)>2$ or $h(\pi)=2$ and the HNN base $H$ has one end then
$H^2(\pi;\mathbb{Z}[\pi])=0$ and so the knot manifold $M(K)$ is aspherical, 
by Theorem 3.5 of \cite{Hi}.
Hence $\pi$ is torsion-free and virtually solvable, 
by Theorem 1.11 of \cite{Hi}.
(Closer examination shows that it must be polycyclic.
See Chapter 16 of \cite{Hi}.)
Otherwise $H$ must have two ends. 
Let $T$ be the maximal finite normal subgroup of $H$. 
Then $\phi(T)=T$, since $\phi$ is injective, and so $T$ is normal in $\pi$.
Hence $T=1$ and $\pi\cong\Phi$, by Theorem 15.2 of \cite{Hi}.

\section{the hantzsche-wendt flat 3-manifold}

The group of affine motions of 3-space is 
$Aff(3)=\mathbb{R}^3\rtimes{GL(3,\mathbb{R})}$.
The action is given by $(v,A)(x)=Ax+v$, for all $x\in\mathbb{R}^3$.
Therefore $(v,A)(w,B)=(v+Aw,AB)$.

Let $\{e_1,e_2,e_3\}$ be the standard basis of $\mathbb{R}^3$,
and let $X,Y,Z\in{GL(3,\mathbb{Z})}$ be the diagonal matrices
$X=diag[1,-1,-1]$, $Y=diag[-1,1,-1]$ and $Z=diag[-1,-1,1]$.
Let $x=(\frac12e_1,X)$, $y=(\frac12(e_2-e_3),Y)$ and 
$z={(\frac12(e_1-e_2+e_3),Z)}$.
The subgroup of $Aff(3)$ generated by $x$ and $y$ is
the Hantzsche-Wendt flat 3-manifold group $G_6$,
with presentation
\[\langle{x,y,z}\mid{xy^2x^{-1}y^2=yx^2y^{-1}x^2=1,~z=xy}\rangle.\]
The translation subgroup $T=G_6\cap\mathbb{R}^3$ is free abelian,
with basis $\{x^2,y^2,z^2\}$.
(This is the maximal abelian normal subgroup of $G_6$.)
The holonomy group $H=\{I,X,Y,Z\}\cong(Z/2Z)^2$
is the image of $G_6$ in $GL(3,\mathbb{R})$.
(Thus $H\cong{G_6/T}$.)
We may clearly take $\{1,x,y,z\}$ as
coset representatives for $H$ in $G_6$.
The commutator subgroup $G_6'$ is free abelian, 
with basis $\{x^4,y^4,x^2y^2z^{-2}\}$.
Thus $2T<G_6'<T$, $T/G_6'\cong(Z/2Z)^2$ and $G_6'/2T\cong{Z/2Z}$.

The orbit space $HW={G_6}\backslash\mathbb{R}^3$ 
is the {\it Hantzsche-Wendt} flat 3-manifold.

\section{the automorphism group of $G_6$}

Every automorphism of the fundamental group of a flat $n$-manifold 
is induced by conjugation in $Aff(n)$,
by a theorem of Bieberbach.
Hence $Aut(G_6)\cong{N/C}$ and $Out(G_6)\cong{N/CG_6}$,
where $C=C_{Aff(3)}(G_6)$ and $N=N_{Aff(3)}(G_6)$.

If $(v,A)\in{Aff(3)}$ commutes with all elements of $G_6$
then $AB=BA $ for all $B\in{H}$, so $A$ is diagonal, and 
$v+Aw=w+Bv$ for all $(w,B)\in{G_6}$.
Taking $B=I$, we see that $Aw=w$ for all $w\in\mathbb{Z}^3$, so $A=I$,
and then $v=Bv$ for all $B\in{H}$, so $v=0$.
Thus $C=1$, and so $Aut(G_6)\cong{N}$.

If $(v,A)\in{N}$ then $A\in{N_{GL(3,\mathbb{R})}(H)}$ 
and $A$ preserves $T=\mathbb{Z}^3$, so $A\in{N_{GL(3,\mathbb{Z})}(H)}$.
Therefore $W=AXA^{-1}$ is in $H$.
Hence $WA=AX$ and so $WAe_1=Ae_1$ is 
up to sign the unique basis vector fixed by $W$.
Applying the same argument to $AYA^{-1}$ and $AZA^{-1}$,
we see that $N_{GL(3,\mathbb{R})}(H)$ is 
the group of ``signed permutation matrices",
generated by the diagonal matrices and permutation matrices.
Let
\[ P=
\left(
\begin{matrix}
0& 1 & 0\\
1 & 0 &0\\
0& 0& -1
\end{matrix}
\right)
\quad \mathrm{and}\quad 
J=
\left(
\begin{matrix}
0& 1 &0\\
0&0&-1\\
1& 0& 0
\end{matrix}
\right).
\]
If $A$ is a diagonal matrix in $GL(3,\mathbb{Z})$ then $(0,A)\in{N}$.
Thus $\tilde{a}=(0,-X)$, $\tilde{b}=(0,-Y)$ and $\tilde{c}=(0,-Z)$ are in $N$.
It is easily seen that $N\cap\mathbb{R}^3=\frac12\mathbb{Z}^3$,
with basis $\tilde{d}=(\frac12{e_1},I)$, $\tilde{e}=(\frac12{e_2},I)$ and 
$\tilde{f}=(\frac12{e_3},I)$.
It is also easily verified that $\tilde{i}=(-\frac14e_3,P)$ 
and $\tilde{j}=(\frac14(e_1-e_2),J)$ are in $N$,
and that $N$ is generated by 
$\{\tilde{a},\tilde{b},\tilde{c},\tilde{d},\tilde{e},
\tilde{f},\tilde{i},\tilde{j}\}$.

The natural action of $N$ on $\mathbb{R}^3$ is isometric,
since ${N_{GL(3,\mathbb{R})}(H)<O(3)}$,
and so $N/G_6$ acts isometrically on the orbit space $HW$.
In fact every isometry of $HW$ lifts to an affine transformation of
$\mathbb{R}^3$ which normalizes $G_6$, and so $Isom(HW)\cong{Out(G_6)}$.
The isometries which preserve the orientation
are represented by pairs $(v,A)$ with $\det(A)=1$. 
(Thus $\tilde{a},\tilde{b},\tilde{c}$ and $\tilde{j}$ 
represent orientation reversing isometries.)

In Chapter 8.\S2 of \cite{Hi} we showed that $Aut(G_6)$ 
is generated by the automorphisms $a,b,c,d,e,f,i$ and $j$ 
which send $x$ to $x^{-1},x,x,x,y^2x$,$z^2x$, $y$, $z$
and $y$ to $y, y^{-1}, z^2y, x^2y,y,z^2y, x,x$, respectively.
These automorphisms are induced by conjugation by 
$\tilde{a},\tilde{b},\tilde{c},\tilde{d},\tilde{e},\tilde{f},\tilde{i}$ 
and $\tilde{j}$, and we shall henceforth drop the tildes.

The subgroup of $Aut(G_6)$ generated by $\{a,b,c,d,e,f\}$ is normal, 
and is a semidirect product $Z^3\rtimes{(Z/2Z)^3}$ with presentation
\[\langle {a,b,c,d,e,f}\mid
{a^2=b^2=c^2=1,}~a,b,c~\text{commute,}~d,e,f~\text{commute,}\]
\[ada=d^{-1}\!,~ae=ea,~af=fa,
~bd=db,~beb=e^{-1}\!,~bf=fb,\]
\[cd=dc,~ce=ec,~cfc=f^{-1}\rangle.\]
This subgroup contains the inner automorphisms $c_x=bcd$,
$c_y=acef$ and $c_z=c_xc_y$ determined by conjugation by $x$ and $y$.
In particular, $c_x^2=d^2$, $c_y^2=e^2$ and $c_z^2=f^2$.
Adjoining the generator $j$ gives another normal subgroup,
in which $j^3=abce$, so $j^6=1$, and $j$ acts on
$\langle {a,b,c,d,e,f}\rangle$ as follows:
\[jaj^{-1}=c,~jbj^{-1}=ad^{-1},~jcj^{-1}=be,~jdj^{-1}=f,\]
\[jej^{-1}=d,~jfj^{-1}=e^{-1}.\]
This subgroup has index 2 in $Aut(G_6)$.
The remaining generator $i$ is an involution ($i^2=1$),
and there are further relations
\[idi=e,~iei=d,~ifi=f^{-1}\!,~iai=b,~ibi=a,~ici=cf,~jiji=d.\]
Therefore $Out(G_6)$ has the presentation
\begin{gather*} 
\langle{a,b,c,e,i,j}\mid{a^2=b^2=c^2=e^2=i^2=j^6=1,}
~a,b,c,e~\text{commute,}\\
iai=b,\medspace{ici=ae,}\medspace{jaj^{-1}=c},\medspace{jbj^{-1}=abc},
\medspace{jcj^{-1}=be,}\medspace{jej^{-1}=bc,}\\
{j^3=abce},\medspace (ji)^2=bc\rangle.
\end{gather*} 
(The images of $d$ and $f$ are represented by $bc$ and $ace$, respectively.)
The natural homomorphism from $Out(G_6)$ to 
$Aut(G_6/T)\cong{GL(2,\mathbb{F}_2)}$ is onto, 
as the images of $i$ and $j$ generate $GL(2,\mathbb{F}_2)$,
and its kernel is the subgroup generated by $\{a,b,c,e\}$.
Thus $Out(G_6)$ has order 96.

\smallskip
\noindent{\sl Comparison with \cite{Zn}}.
As the group element labeled $z$ and the automorphisms labeled $a,b,c$ by
Zimmermann differ from ours, 
we shall add the subscript ``$Z$" for clarity.
The presentation for $G_6$ used in \cite{Zn} reduces to
\[\langle{x,y,z_Z}\mid{xy^2x^{-1}y^2=yx^2y^{-1}x^2=1,~z_Zyx=x^2z_Z^2=z_Z^2x^2}
\rangle.\]
Thus $z_Z=yx^{-1}$, so $z_Z=y^2z^{-1}$, and $z_Z^2=z^{-2}$.
His choice of representatives for a generating set for $Out(G_6)$ 
is $\{a_Z,b_Z,c_Z,I,S,T\}$,
where $a_Z=d^{-1}$, $b_Z=e^{-1}$, $c_Z=f^{-1}$, 
$I=xade$, $S=ideab$ and $T=j^{-1}i^{-1}xadf$.
He observes also that $Out(G_6)$ is an extension of $S_3\times{Z/2Z}$
by the normal subgroup $(Z/2Z)^3$ generated by $\{d,e,f\}$,
but the extension does not split, since the centre of $Out(G_6)$
is generated by the image of the involution $ab$, and thus is too small.

\section{2-knots with $\pi'\cong{G_6}$ are not reflexive}

Since $G_6$ is solvable and $H_1(G_6)\cong(Z/4Z)^2$,
an automorphism $(v,A)$ of $G_6$ is meridianal 
if and only if its image in $Aut(G_6/T)\cong{GL(2,2)}$ has order 3.
Thus its image in $Out(G_6)$ is conjugate to $[j]$, $[j]^{-1}$, 
$[ja]$ or $[jb]$.
The latter pair are orientation-preserving 
and each is conjugate to its inverse (via $[i]$).
However $(ja)^3=1$ while $(jb)^3=de^{-1}f$, so $[jb]^3=[ab]\not=1$.
Thus $[ja]$ is not conjugate to $[jb]^\pm$, and the knot groups
$G(+)=G_6\rtimes_{[ja]}Z$ and $G(-)=G_6\rtimes_{[jb]}Z$ are distinct.
The corresponding knot manifolds are the mapping tori 
of the isometries of $HW$ determined by $[ja]$ and $[jb]$, 
and are flat 4-manifolds.

We may use the geometry of flat manifolds to adapt 
the argument of Lemma 18.3 of \cite{Hi} 
to our present situation, as follows.

\begin{theorem}
Let $K$ be a $2$-knot with group $G(+)$ or $G(-)$.
Then $K$ is not reflexive.
\end{theorem}

\begin{proof}
The knot manifold $M=M(K)$ is a flat 4-manifold,
$M\cong\mathbb{R}^4/\pi$ say, 
by Theorem 8.1 and the discussion in Chapter 16 of \cite{Hi}.
The weight orbit of $K$ may be represented by a geodesic 
simple closed curve $C$ through the basepoint $P$ of $\mathbb{R}^4/\pi$.
Let $\gamma$ be the image of $C$ in $\pi$.

Let $h$ be a self-homeomorphism of $\mathbb{R}^4/\pi$ which fixes $C$
pointwise.
Since $M$ is aspherical $h$ is based-homotopic to an affine diffeomorphism 
$\alpha$, and then $\alpha_*(\gamma)=h_*(\gamma)=\gamma$.
Let $\widehat{M}\cong\mathbb{R}^3\times{S^1}$ be the covering space 
corresponding to the subgroup $\langle\gamma\rangle\cong{Z}$,
and fix a lift $\widehat{C}$.
A homotopy from $h$ to $\alpha$ lifts to a proper homotopy
between the lifts $\widehat{h}$ and $\widehat{\alpha}$ to 
self-homeomorphisms fixing $\widehat{C}$.
Now the behaviour at $\infty$ of these maps is determined by
the behaviour near the fixed point sets, 
as in \cite{CS,HP} or Lemma 18.3 of \cite{Hi}.
Since the affine diffeomorphism $\widehat{\alpha}$ does not change
the framing of the normal to $\widehat{C}$ it follows that
$\widehat{h}$ and $h$ do not change the normal framings either.
\end{proof}

\section{symmetries of 2-knots with group $G(+)$}

The orthogonal matrix $-JX$ is a rotation though $\frac{2\pi}3$
about the axis in the direction $e_1+e_2-e_3$.
The fixed point set of the isometry $[ja]$ of $HW$ is
the image of the line $\lambda(s)=s(e_1+e_2-e_3)-\frac14e_2$.
The knots corresponding to the canonical section are the 
3-twist spin of the figure eight knot $\tau_34_1$ 
and its Gluck reconstruction $\tau_34_1^*$.
The knot $\tau_34_1$ is $\pm$amphicheiral and invertible \cite{Li}.
We shall show that $\tau_34_1$ is {\it strongly\/} $\pm$amphicheiral, 
but not strongly invertible.
We shall also show that none of the other 2-knots with
group $G(+)$ are amphicheiral or invertible.

\begin{theorem}
Let $\pi=G(+)$.
Then every strict weight orbit representing a given generator $t$ 
for $\pi/\pi'$ contains an unique element of the form $x^{2n}t$.
\end{theorem}

\begin{proof}
If $t\in\pi$ represents a generator of $\pi/\pi'\cong{Z}$
it is a weight element, since $\pi$ is solvable.
Suppose that $c_t=ja$.
If $u$ is another weight element with $[c_u]=[ja]$ 
then $c_u$ is conjugate in $Aut(G_6)$
to $c_{gt}$, for some $g\in\pi''=G_6'$,
by Theorem 14.1 of \cite{Hi}.
Suppose that $g=x^{2m}y^{2n}z^{2p}$.
Let $\lambda(g)=m+n-p$ and $w=x^{2n}y^{2p}$.
Then $w^{-1}gtw=x^{2\lambda(g)}t$.
On the other hand, if $\psi\in{Aut(G_6)}$ then 
$\psi{c_{gt}}\psi^{-1}=c_{ht}$ for some $h\in{G_6'}$ 
if and only if the images of $\psi$ and $ja$ in $Aut(G_6)/G_6'$ commute.
If so, $\psi$ is in the subgroup generated by $\{def^{-1},jb,ce\}$.
It is easily verified that $\lambda(h)=\lambda(g)$
for any such $\psi$. 
(It suffices to check this for the generators.)

Thus $x^{2n}t$ is a weight element representing $[ja]$,
for all $n\in\mathbb{Z}$, and $x^{2m}t$ and $x^{2n}t$ are in 
the same strict weight orbit if and only if $m=n$.
\end{proof}

Note that $\lambda$ is given by dot product with the axis of $-JX$.

\begin{lemma}
If $n=0$ then $C_{Aut(G_6)}(ja)$ and $N_{Aut(G_6)}(\langle{ja}\rangle)$
are generated by $\{ja,def^{-1},abce\}$ and $\{ja,ice,abce\}$,
respectively.
The subgroup which preserves the orientation of $\mathbb{R}^3$ 
is generated by $\{ja,ice\}$.

If $n\not=0$ then $N_{Aut(G_6)}(\langle{d^{2n}ja}\rangle)=
C_{Aut(G_6)}(d^{2n}ja)$ and is
generated by $\{d^{2n}ja,def^{-1}\}$.
This subgroup acts orientably on $\mathbb{R}^3$.
\end{lemma}

\begin{proof}
This is straightforward.
(Note that $abce=j^3$ and $def^{-1}=(ice)^2$.)
\end{proof}

\begin{lemma}
The mapping torus $M([ja])$ has an orientation reversing involution 
which fixes a canonical section pointwise,
and an orientation reversing involution which fixes a
canonical section setwise but reverses its orientation.
There is no orientation preserving involution of $M$ which reverses 
the orientation of any section.
\end{lemma}

\begin{proof}
Let $\omega=abcd^{-1}f=abce(ice)^{-2}$\!,
and let $p=\lambda(\frac14)=\frac14(e_1-e_3)$.
Then $\omega=(2p,-I_3)$, $\omega^2=1$,
$\omega{ja}=ja\omega$ and $\omega(p)=ja(p)=p$.
Hence $\Omega=m([\omega])$ is an orientation reversing involution 
of $M([ja])$ which fixes the canonical section determined by 
the image of $p$ in $HW$.

Let $\Psi([f,s])=[[iab](f),1-s]$ for all $[f,s]\in{M([ja])}$.
This is well-defined, since $(iab){ja}(iab)^{-1}=(ja)^{-1}$, 
and is an involution, since $(iab)^2=1$.
It is clearly orientation reversing, 
and since $iab(\lambda(\frac18))=\lambda(\frac18)$ 
it reverses the section determined by the image 
of $\lambda(\frac18)$ in $HW$.

On the other hand, $\langle{ja},ice\rangle\cong{Z/3Z}\rtimes_{-1}Z$,
and the elements of finite order in this group do not invert $ja$.
\end{proof}

\begin{theorem}
Let $K$ be a $2$-knot with group $G(+)$ and weight element $u=x^{2n}t$,
where $t$ is the canonical section.
If $n=0$ then $K$ is strongly $\pm$amphicheiral,
but is not strongly invertible.
If $n\not=0$ then $K$ is neither amphicheiral nor invertible.
\end{theorem}

\begin{proof}
Suppose first that $n=0$.
Since $-JX$ has order 3 it is conjugate in $GL(3,\mathbb{R})$
to a block diagonal matrix $\Lambda(-JX)\Lambda^{-1}=\left(\smallmatrix 1&0\\
0&R(\frac{2\pi}3)\endsmallmatrix\right)$,
where $R(\theta)\in{GL(2,\mathbb{R})}$ is rotation through $\theta$.
Let $R_s=R(\frac{2\pi}3s)$ and 
$\xi(s)=((I_3-A_s)p,A_s)$, where 
$A_s= \Lambda^{-1}\left(\smallmatrix 1&0\\
0&R_s\endsmallmatrix\right)\Lambda$, for $s\in\mathbb{R}$.
Then $\xi$ is a 1-parameter subgroup of $Aff(3)$,
such that $\xi(s)(p)=p$ and $\xi(s)\omega=\omega\xi(s)$ for all $s$.
In particular, $\xi|_{[0,1]}$ is a path from $\xi(0)=1$ 
to $\xi(1)=ja$ in $Aff(3)$.
Let $\Xi:\mathbb{R}^3\times{S^1}\to{M(ja)}$ be the homeomorphism given by
$\Xi(v,e^{2\pi{is}})=[\xi(v),s]$ for all $(v,s)\in \mathbb{R}^3\times[0,1]$.
Then $\Xi^{-1}\Omega\Xi=\omega\times{id_{S^1}}$ and so
$\Omega$ does not change the framing.
Therefore $K$ is strongly $-$amphicheiral.
 
Similarly, if we let $\zeta(s)=((I_3-A_s)\lambda(\frac18),A_s)$
then $\zeta(s)(\lambda(\frac18))=\lambda(\frac18)$ and
$iab\zeta(s)iab=\zeta(s)^{-1}$, for all $s\in\mathbb{R}$,
and $\zeta|_{[0,1]}$ is a path from 1 to $ja$ in $Aff(3)$.
Let $Z:\mathbb{R}^3\times{S^1}\to{M(ja)}$ be the homeomorphism given by
$Z(v,e^{2\pi{is}})=[\zeta(s)(v),s]$ for all $(v,s)\in \mathbb{R}^3\times{S^1}$.
Then 
\[Z^{-1}\Psi{Z}(v,z)=(\zeta(1-s)^{-1}iab\zeta(s)(v),z^{-1})\]
\[=(\zeta(1-s)^{-1}\zeta(s)^{-1}iab(v),z^{-1})=
((ja)^{-1}iab(v),z^{-1})\]
for all $(v,z)\in \mathbb{R}^3\times{S^1}$.
Hence $\Psi$ does not change the framing, 
and so $K$ is strongly $+$amphicheiral.
However it is not strongly invertible, by Lemma 6.

If $n\not=0$ every such self-homeomorphism $h$ 
preserves the orientation and fixes the meridian, by Lemma 5,
and so $K$ is neither amphicheiral nor invertible.
\end{proof}

\section{symmetries of 2-knots with group $G(-)$}

A similar analysis applies when the knot group is $G(-)$,
i.e., when the meridianal automorphism is
$jb=(\frac14(e_1-e_2),-JY)$.
(The orthogonal matrix $-JY$ is now a rotation though $\frac{2\pi}3$
about the axis in the direction ${e_1-e_2+e_3}$.)
All 2-knots with group $G(-)$ are fibred,
and the characteristic map $[jb]$ has finite order, 
but none of these knots are twist-spins, as we shall show below.

\begin{theorem}
Let $\pi=G(-)$.
Then every strict weight orbit representing a given generator $t$ 
for $\pi/\pi'$ contains an unique element of the form $x^{2n}t$.
\end{theorem}

\begin{proof}
The proof is very similar to that of Theorem 4.
The main change is that we should define the homomorphism $\lambda$
by $\lambda(x^{2m}y^{2n}z^{2p})=m-n+p$. 
\end{proof}

\begin{cor}
No $2$-knot with group $G(-)$ is a twist-spin.
\end{cor}

\begin{proof}
Suppose that $G(-)$ is the group of the $r$-twist-spin
of a classical knot.
Then the $r$th power of a meridian is central.
The power $(x^{2n}t)^r$ is central in $G(-)$ if and only if $(d^{2n}jb)^r=1$
in $Aut(G_6)$.
But $(d^{2n}jb)^3=d^{2n}f^{2n}e^{-2n}(jb)^3=(de^{-1}f)^{2n+1}$.
Therefore $d^{2n}jb$ has infinite order,
and so $G(-)$ is not the group of a twist-spin.
\end{proof}

\begin{lemma}
If $n=0$ then $C_{Aut(G_6)}(jb)$ and $N_{Aut(G_6)}(\langle{jb}\rangle)$
are generated by $\{jb\}$ and $\{jb,i\}$,
respectively.
If $n\not=0$ then 
$N_{Aut(G_6)}(\langle{d^{2n}jb}\rangle)=C_{Aut(G_6)}(d^{2n}jb)$ 
and is generated by $\{d^{2n}jb,de^{-1}f\}$.
These subgroups act orientably on $\mathbb{R}^3$.
\qed
\end{lemma}

The isometry $[jb]$ has no fixed points in $G_6\backslash\mathbb{R}^3$.
We shall defined a preferred section as follows.
Let $\gamma(s)=\frac{2s-1}8(e_1-e_2)-\frac18e_3$, for $s\in\mathbb{R}$.
Then $\gamma(1)=jb(\gamma(0))$, 
and so $\gamma|_{[0,1]}$ defines a section of $p_{[jb]}$.
We shall let the image of $(\gamma(0),0)$ be the basepoint for $M([jb])$.

\begin{theorem}
Let $K$ be a $2$-knot with group $G(-)$ and weight element $u=x^{2n}t$,
where $t$ is the canonical section.
If $n=0$ then $K$ is strongly $+$amphicheiral but not invertible.
If $n\not=0$ then $K$ is neither amphicheiral nor invertible.
\end{theorem}

\begin{proof}
Suppose first that $n=0$.
Since $i(\gamma(s))=\gamma(1-s)$ for all $s\in\mathbb{R}$
the section defined by $\gamma|_{[0,1]}$ is fixed setwise 
and reversed by the orientation reversing involution 
${[f,s]\mapsto[[i](f),1-s]}$.
Let $B_s$ be a 1-parameter subgroup of $O(3)$ such that $B_1=jb$.
Then we may define a path from 1 to $jb$ in $Aff(3)$ by setting
$\zeta(s)={((I_3-B_s)\gamma(s),B_s)}$ for $s\in\mathbb{R}$.
We see that $\zeta(0)=1$, $\zeta(1)=jb$,
$\zeta(s)(\gamma(s))=\gamma(s)$ 
and $i\zeta(s)i=\zeta(s)^{-1}$ for all $0\leq{s}\leq1$.
As in Theorem 5 it follows that the involution does not change the framing
and so $K$ is strongly $+$amphicheiral.

The other assertions follow from Lemma 9, as in Theorem 7.
\end{proof}

In particular, only $\tau_34_1, \tau_34_1^*$ and the knots obtained by
surgery on the section of $M([jb])$ defined by $\gamma|_{[0,1]}$
admit orientation-changing symmetries.

\section{normal forms for meridianal automorphisms of $\Gamma(e,\eta)$}

Let $M(e,\eta)$ be the 2-fold branched covering of $S^3$, 
branched over a Montesinos knot 
$k(e,\eta)=K(0|e;(3,\eta),(3,1),(3,1))$,
with $e$ even and $\eta=\pm1$.
This 3-manifold is Seifert fibred over 
the flat 2-orbifold $S(3,3,3)$,
and $\Gamma(e,\eta)=\pi_1(M(e,\eta))$ has a presentation
\begin{equation*}
\langle 
h,x,y,z\mid x^3=y^3=z^{3\eta}=h,\medspace xyz=h^e\rangle,
\end{equation*}
for some $\eta=\pm1$.
Let $u=z^{-1}x$, $v=xz^{-1}$ and $q=3e-\eta-2$. 
Then $\Gamma(e,\eta)$ also has the presentation
\[\langle{u,v,z}\mid{zuz^{-1}=v},
~zvz^{-1}=v^{-1}u^{-1}z^{3\eta-3},
~vuv^{-1}u^{-1}=z^{3\eta{q}}\rangle.\]
The image of $z^{3\eta}$ in $\Gamma(e,\eta)$ generates the centre 
$\zeta\Gamma(e,\eta)$,
and $P=\Gamma(e,\eta)/\zeta\Gamma(e,\eta)$ 
is the orbifold fundamental group of $S(3,3,3)$.

An automorphism $\phi$ of $\Gamma(e,\eta)$ must preserve 
characteristic subgroups such as the centre $\zeta\Gamma(e,\eta)$ 
(generated by $z^3$) and the maximal nilpotent normal subgroup 
$\sqrt{\Gamma(e,\eta)}$ (generated by $u$, $v$ and $z^3$).
Let $F$ be the subgroup of $Aut(\Gamma(e,\eta))$ 
consisting of automorphisms which induce the identity on 
$\Gamma(e,\eta)/\sqrt{\Gamma(e,\eta)}\cong{Z/3Z}$ and 
$\sqrt{\Gamma(e,\eta)}/\zeta\Gamma(e,\eta)\cong{Z}^2$.
Automorphisms in $F$ also fix the centre, and are of the form $k_{m,n}$,
where
\[{k_{m,n}(u)=uz^{3\eta{s}}},~{k_{m,n}(v)=vz^{3\eta{t}}}\quad
\mathrm{and}\quad{k_{m,n}(z)=z^{3\eta{p}+1}u^mv^n},\]
for $(m,n)\in\mathbb{Z}^2$. 
These formulae define an automorphism if and only if
\[s-t=-nq,\quad {s+2t}=mq\quad\mathrm{and}\quad 
6p=(m+n)((m+n-1)q+2(\eta-1)).\]
In particular, conjugation by $u$ and $v$ give
$c_u=k_{-2,-1}$ and $c_v=k_{1,-1}$, respectively.
If $\eta=1$ then $q=3e$, so $s=(m-2n)e$, $t=-(m+n)e$
and $p=\binom{m+n}2e$ are integers for all $m,n\in\mathbb{Z}$.
In this case $F\cong\mathbb{Z}^2$ is generated by $k=k_{1,0}$ and $c_u$.
If $\eta=-1$ then $m+n\equiv0$ {\it mod} $(3)$,
and $F$ is generated by $c_u$ and $c_v$.
In this case $F$ has index 3 in $\mathbb{Z}^2$.

We may define automorphisms $b$ and $r$ by the formulae:
\[{b(u)=v^{-1}z^{3\eta{e}-3}},\quad 
b(v)=uvz^{3\eta(e-1)}\quad\mathrm{and}\quad {b(z)=z};\quad\mathrm{and}\]
\[{r(u)=v^{-1}},\quad{r(v)=u^{-1}}\quad\mathrm{and}\quad{r(z)=z^{-1}}.\]
It is easily checked that ${b^6=r^2=(br)^2=1}$ and that
conjugation by $z$ gives $c_z=b^4$.
Since $\Gamma(e,\eta)/\Gamma(e,\eta)'$ is finite, 
$Hom(\Gamma(e,\eta),\zeta\Gamma(e,\eta))=0$,
and so the natural homomorphism from $Aut(\Gamma(e,\eta))$ 
to $Aut(P)$ is injective.
If $\eta=+1$ this homomorphism is an isomorphism,
and $Aut(\Gamma(e,1))=$ has a presentation
\[\langle{b,c_u,k,r}\mid
{b^6=r^2=(br)^2=1},~c_uk=kc_u,~bc_ub^{-1}=c_u^{-1}k^{-3}\!,\]
\[~bkb^{-1}=c_uk^2,~rc_ur=c_uk^3,~rkr=k^{-1}\rangle.\]
(Here $c_v=c_uk^3$.) 
On the other hand, $Aut(\Gamma(e,-1))$ has a presentation
\[\langle{b,c_u,c_v,r}\mid
{b^6=r^2=(br)^2=1},~c_uc_v=c_vc_u,~bc_ub^{-1}=c_v^{-1}\!,\]
\[~bc_vb^{-1}=c_uc_v,~rc_ur=c_v,~rc_vr=c_u\rangle.\]
Hence $Out(\Gamma(e,1))\cong{S_3}\times{Z/2Z}$,
while $Out(\Gamma(e,-1))\cong(Z/2Z)^2$.

In each case an automorphism $\phi$ is meridianal if and only if
$[\phi]$ is conjugate to $[r]$, and so there is an unique corresponding knot
group $\pi(e,\eta)=\Gamma(e,\eta)\rtimes_r{Z}$.

\section{embeddings in the affine group}

The group $P$ embeds as a discrete subgroup of $Isom(\mathbb{E}^2)$,
via $u\mapsto(e_1,I_2)$, $v\mapsto(e_2,I_2)$ and $z\mapsto(0,-\beta)$,
where $\beta=\left(\smallmatrix 0&1\\
-1&1\endsmallmatrix\right)$.
The images of $u$ and $v$ in $P$ 
form a basis for the translation subgroup $T(P)\cong Z^2$,
and $P=T(P)\rtimes_{-\beta}(Z/3Z)$.
It is easily seen that $C_{Aff(2)}(P)=1$,
and so $Aut(P)\cong{N_{Aff(2)}}(P)$.
If $(v,A)\in{N}_{Aff(2)}(P)$ then $(I_2+\beta)v\in\mathbb{Z}^2$ 
and $A$ is in the subgroup $D$ of $GL(2,\mathbb{R})$
generated by the matrices $\beta$ and 
$\rho=\left(\smallmatrix 0&-1\\
-1&0\endsmallmatrix\right)$, 
which has order 12.
Thus $Aut(P)=(I_2+\beta^{-1})T(P)\rtimes{D}$.
Hence $Out(P)\cong{D}\cong{S_3}\times{Z/2Z}$,
where the first factor is generated by the images of $u$ and $(0,-I_2)$
and the second factor is generated by the image of $(0,-\rho)$.

Let $Nil<GL(3,\mathbb{R})$ be the group of $3\times3$ upper triangular
matrices 
\[[x,y,w]=
\left(
\begin{matrix}
1& x& w\\
0&1& y\\
0&0&1
\end{matrix}
\right),
\]
and let $Aut(Nil)$ be the group of Lie automorphisms.
As a set, $Aut(Nil)$ is the cartesian product
$GL(2,\mathbb{R})\times\mathbb{R}^2$, with 
$(A,\mu)= (\left(\smallmatrix a& c\\
          b& d\endsmallmatrix\right),(\mu_1,\mu_2))$ 
acting via $(A,\mu)([x,y,w])=$
\begin{equation*}
[ax+cy,bx+dy,\mu_1x+\mu_2y+(ad-bc)w+bcxy+\frac{ab}2x(x-1)+\frac{cd}2y(y-1)].
\end{equation*}
All such automorphisms are orientation preserving.
The product of $(A,\mu)$ with 
$(B,\nu)=(\left(\smallmatrix g&j\\
h&k\endsmallmatrix\right),(n_1,n_2))$
is \[(A,\mu)\circ(B,\nu)\! =\! (AB,\mu B+det(A)\nu+\frac12\eta(A,B)),\]
where
\[\eta(A,B)=(abg(1-g)+cdh(1-h)-2bcgh,abj(1-j)+cdk(1-k)-2bcjk).\]

Let $Aff(Nil)=Nil\rtimes{Aut(Nil)}$.
Then $Aff(Nil)$ acts on the open 3-manifold $Nil\cong\mathbb{R}^3$
by $(n,\sigma)(n')=n\sigma(n')$.
The abelianization $Nil\to\mathbb{R}^2=Nil/\zeta{Nil}$ 
extends to an epimorphism $p:Aff(Nil)\to{Aff(2)}$,
given by 
$p(n,A,\mu)=(\left(\smallmatrix x\\
y\endsmallmatrix\right),A)$ 
for $n=[x,y,w]\in{Nil}$,
$A\in{GL(2,\mathbb{R})}$ and $\mu\in\mathbb{R}^2$.
We may embed $\Gamma(e,\eta)$ in $Aff(Nil)$ by
\[u\mapsto([1,0,0],\iota),\quad{v}\mapsto([0,1,0],\iota)\quad\mathrm{and}\quad
z\mapsto([0,0,\frac{-1}{3q}],\alpha),\] 
where $\iota=id_{Nil}$ and $\alpha=(-\beta,(0,\frac{\eta-1}q))$.
(Note that $vuv^{-1}u^{-1}\mapsto([0,0,-1],\iota)$.)
Let $N=N_{Aff(Nil)}(\Gamma(e,\eta))$ and $C=C_{Aff(Nil)}(\Gamma(e,\eta))$.
As in the flat case, 
$Aut(\Gamma(e,\eta))\cong{N/C}$ and 
$Out(\Gamma(e,\eta))\cong{N/C}\Gamma(e,\eta)$.

It is easily seen that $C=\zeta{Nil}=\{([0,0,z],\iota)\mid{z\in\mathbb{R}}\}$.
If $n=[x,y,w]$ and $(n,A,\mu)\in{N}$ then 
$(\left(\smallmatrix x\\
y\endsmallmatrix\right),A)\in{N_{Aff(2)}(P)}$,
so $A\in{D}$ and
$\left(\smallmatrix x\\
y\endsmallmatrix\right)\in(I_2+\beta)^{-1}\mathbb{Z}^2$.
If $A=I_2$ then $(n,I_2,\mu)$ is in $N$ if and only if
it normalizes $\sqrt{\Gamma(e,\eta)}$
and $(n,I_2,\mu)z=(n',\iota)z(n,I_2,\mu)$
for some $n'\in\sqrt{\Gamma(e,\eta)}$.
The latter condition implies that $(I_2,\mu)\alpha=\alpha(I_2,\mu)$,
and so $\mu(\beta+I_2)=0$.
Thus we must have $\mu=0$ and $(I_2,\mu)=\iota$.
The remaining conditions then imply that $x,y\in\frac1q\mathbb{Z}$.
If $\eta=1$ (so $q=3e$) this satisfied by all
$\left(\smallmatrix x\\
y\endsmallmatrix\right)\in(I_2+\beta)^{-1}\mathbb{Z}^2<\frac13\mathbb{Z}^2$.
If $\eta=-1$ then $x,y\in\mathbb{Z}$.
Thus the natural map from $Aut(\Gamma(e,\eta))$ to $Aut(P)$ 
is an isomorphism if $\eta=1$, and has image of index 3 if $\eta=-1$.

\section{2-knots with group $\pi(e,\eta)$}

Let $R=([0,0,0],\rho,(0,0))$ in $Aff(Nil)$.
Then $R^2=1$ and $R([x,y,z])$ $=[-y,-x,-z]$ for all $[x,y,z]\in{Nil}$.
The fixed point set of the action of $R$ on $Nil$ is the
connected curve $\{[s,-s,0]\mid{s\in\mathbb{R}}\}$.
Thus the fixed point set of the involution $[R]$ of $M(e,\eta)$ 
induced by $R$ is connected and nonempty.
The corresponding 2-knot is $\tau_2k(e,\eta)$.
This is reflexive and $+$amphicheiral \cite{Li}.

\begin{theorem}
The knot $K=\tau_2k(e,\eta)$ is reflexive and strongly 

\noindent$+$amphicheiral, but is not invertible.
\end{theorem}

\begin{proof}
Let $S([m,s])=[b^3(m),s]$ and $h([m,s])=[m,1-s]$ for $m\in{M(e,\eta)}$
and $0\leq{s}\leq1$.
Then $S$ and $h$ define commuting involutions of $M([R])$,
which each fix the canonical section setwise.

As remarked in \S3, in order to determine how these involutions
affect the framing we may pass to the 
irregular covering space $M(R)=$ ${Nil\times_RS^1}$.
We shall identify the space $Nil$ with $\mathbb{R}^3$, 
in the obvious way.

Let $R(\theta)\in{GL(2,\mathbb{R})}$ be rotation through $\theta$,
and let $P=\left(\smallmatrix R(\frac\pi4) &0\\
0&1\endsmallmatrix\right)\in{GL(3,\mathbb{R})}$.
Then $PRP^{-1}=diag[1,-1,-1]$.
We may isotope $PRP^{-1}$ back to the identity,
via $Q_s=\left(\smallmatrix 1&0\\
0& R(s\pi)\endsmallmatrix\right)$, for $0\leq{s}\leq1$.
Let $Q:\mathbb{R}^3\times{S^1}\to{M(PRP^{-1})}$
be the homeomorphism given by $Q(v,e^{2\pi{is}})=[Q_s(v),s]$ 
for all $(v,s)\in\mathbb{R}^3\times[0,1]$.
Then $Q^{-1}hQ((v,z)=(Q_{2s-1}(v),z^{-1})$
for all $(v,z)\in\mathbb{R}^3\times{S^1}$.
After reversing the $S^1$ factor this is just the twist,
and so $h$ changes the framing.
Thus $K$ is reflexive.

The automorphism $b^3$ acts linearly, 
via $b^3([x,y,z])=[-x,-y,z+(e\eta-1)(x+y)]$, 
and so $Pb^3P^{-1}=
\left(\smallmatrix -I_2&0\\
\mu&1\endsmallmatrix\right)$, where $\mu=(e\eta-1,e\eta-1)R(-\frac\pi4)$.
We may isotope $Pb^3P^{-1}$ to $d=diag[-1,-1,1]$ through 
invertible matrices which commute with $PRP^{-1}$. 
Let $D([v,s])=[d(v),s]$. 
Then $S$ and $D$ twist the framing in the same way.
Since $Q^{-1}DQ(v,e^{2\pi{is}})$ 
$=(Q_{-s}dQ_s(v),e^{2\pi{is}})=(dQ_{2s}(v),e^{2\pi{is}})$,
for all $(v,s)\in\mathbb{R}^3\times[0,1]$,
it follows that $S$ changes the framing.

The composite $Sh$ is an involution which reverses the orientation 
and the meridian, but does not twist the framing.
Hence $K$ is strongly $+$amphicheiral.

Since automorphisms of $\Gamma(e,\eta)$ are orientation 
preserving $K$ is not $-$amphicheiral or invertible.
\end{proof}

The other knots with such groups are less symmetrical.

\begin{theorem}
Let $K$ be a 2-knot with $\pi=\pi{K}\cong\pi(e,\eta)$.
Then $K$ is reflexive, 
and every strict weight orbit representing the canonical section $t$ 
for $\pi/\pi'$ contains an unique element of the form $u^nt$.
If $n\not=0$ then 
$N_{Aut(\Gamma(e,\eta))}(\langle{u^nr}\rangle)=C_{Aut(\Gamma(e,\eta))}(u^nr)=
\langle{u^nr,uv^{-1}}\rangle$, and $u^nr$ is not conjugate to its inverse.
Hence $K$ is neither amphicheiral nor invertible.
\end{theorem}

\begin{proof}
The first assertion follows from Lemma 1.

Recall that $F$ is the group of automorphisms of $\Gamma(e,\eta)$
which induce the identity on $\Gamma(e,\eta)/\sqrt{\Gamma(e,\eta)}$ and 
$\sqrt{\Gamma(e,\eta)}/\zeta\Gamma(e,\eta)$.
If $\psi{r}\psi^{-1}=rk$ for some $k\in{F}$ then we may assume that
$\psi\in{F}$, and then $k\in(I-\rho)F$.
The paremetrization of the weight orbits follows by the argument of Theorem 4, 
with minor changes.
(Note that $t\mapsto{th}$ defines an automorphism of $\pi$.)

If $n\not=0$ then $u^nt$ is not conjugate to its inverse
Hence $K$ is not $+$amphicheiral, 
and since automorphisms of $\Gamma(e,\eta)$ are orientation 
preserving $K$ is neither $-$amphicheiral nor invertible.
\end{proof}

\section{double null concordance}

A $2$-knot $K$ is doubly null-concordant
(or doubly slice) if $K=S^4\cap{U}$,
where $S^4$ is embedded as the equator of $S^5$ 
and $U$ is a trivial 3-knot.
An equivalent condition is that $M(K)$ embeds in $S^1\times{S^{n+2}}$,
via a map which induces the abelianization on $\pi{K}$.
Thus whether $K$ is doubly null concordant depends only on 
(the $h$-cobordism class of) $M(K)$.

The knot $k(0,-1)=9_{46}$ is doubly null-concordant \cite{HK,Su},
and hence so are all of its twist spins.
Since $M(\tau_2k(0,-1))$ is determined up to homeomorphism by its group,
every knot with group $\pi(0,-1)$ must be doubly null-concordant.

The Farber-Levine pairing of a doubly null-concordant knot is hyperbolic,
and so the torsion submodule of $\pi'/\pi''$ is a direct double. 
Moreover, if $F$ is a field then $H^*(M(K)';{F})$ splits 
as a sum of graded subalgebras ${A^*}\oplus{B^*}$, 
with $A^i\cup{A^{3-i}}=B^i\cup{B^{3-i}}=0$ for $i=1$ or 2,
while duality gives a perfect pairing of $A^i$ with $B^{3-i}$ \cite{Le83}.

The Alexander polynomials of Fox's knot,
the Cappell-Shaneson 2-knots and the knots with 
$\pi'\cong\Gamma_q$ for some odd $q\geq1$ are all irreducible
and nonconstant. 
For such knots $H^*(M(K)';\mathbb{Q})$ cannot split as above.
If $\pi\cong{G(\pm)}$ then $\pi'/\pi''\cong(Z/4Z)^2$ is finite.
However it is cyclic as a $\Lambda$-module,
and so the cohomology ring with coefficients $\mathbb{F}_2$
does not split.
If $\pi=\pi(e,\eta)$ and $|q|=|3e-\eta-2|\not=1$ 
then $\pi'/\pi''\cong{Z/3qZ}\oplus{Z/3Z}$ is not a direct double, 
and so the Farber-Levine pairing is not hyperbolic.
Thus none of these knots are doubly null-concordant.

Similarly, knots whose groups have finite commutator subgroup
other than $I^*$ are not doubly slice, since their
Farber-Levine pairings are not hyperbolic.
All 2-knots with perfect commutator subgroup are topologically doubly slice.
However $\tau_53_1$ is not {\it smoothly} doubly slice \cite{Ru}.
(Crude estimates based on the surgery exact sequence show that there
are infinitely many such knots.)

The classification of these knots up to {\it stable\/} 
DNC equivalence is much more subtle, 
and we shall not attempt this here.

All fibred 2-knots may be realized as smooth knots in
a smooth homotopy 4-sphere.
Twist spins (such as $\tau_34_1$ and $\tau_2k(e,\eta)$) are
smooth knots in the standard smooth structure on $S^4$.
In \cite{Gf} is is shown that this holds also for many Cappell-Shaneson 2-knots.
Is this so for the other knots with torsion-free solvable groups?

\end{document}